\documentclass[10pt, reqno]{amsart}
\usepackage{amsmath, amsthm, amscd, amsfonts, amssymb, graphicx, color}
\usepackage[bookmarksnumbered, colorlinks, plainpages]{hyperref}
\hypersetup{colorlinks=true,linkcolor=red, anchorcolor=green, citecolor=cyan, urlcolor=red, filecolor=magenta, pdftoolbar=true}
\usepackage{mathrsfs}

\newtheorem{theorem}{Theorem}[section]
\newtheorem{lemma}[theorem]{Lemma}
\newtheorem{proposition}[theorem]{Proposition}

\theoremstyle{definition}
\newtheorem{definition}[theorem]{Definition}
\newtheorem{example}[theorem]{Example}

\theoremstyle{remark}
\newtheorem{remark}[theorem]{Remark}
\numberwithin{equation}{section}

\begin{document}
\setcounter{page}{1}

\title[Serre $C^*$-algebras]{Remark on Serre $C^*$-algebras}

\author[Nikolaev]
{Igor V. Nikolaev$^1$}

\address{$^{1}$ Department of Mathematics and Computer Science, St.~John's University, 8000 Utopia Parkway,  
New York,  NY 11439, United States.}
\email{\textcolor[rgb]{0.00,0.00,0.84}{igor.v.nikolaev@gmail.com}}


\subjclass[2010]{Primary 14A22; Secondary 46L85.}

\keywords{twisted homogeneous coordinate ring,  crossed product,  Takai duality.}


\begin{abstract}
We study non-commutative algebraic geometry of
  Artin, Serre  and Tate in terms of the operator algebras.  
Namely, we define the Serre  $C^*$-algebra $\mathcal{A}_X$
of a projective variety $X$ as the norm-closure of a  representation of 
the  twisted homogeneous coordinate ring of $X$ 
by   the linear operators  on a Hilbert   space $\mathcal{H}$. 
It is proved that $X$ is homeomorphic   to the space of all   
irreducible representations  of   the crossed product of $\mathcal{A}_X$  by 
an  automorphism  of $\mathcal{A}_X$. 
The case of rational elliptic curves $X$ is considered in detail.  
\end{abstract}

\maketitle

\section{Introduction}
Since  the  rings of polynomials  are commutative, it is recognized  that
 algebraic geometry must be based on  commutative 
algebra.    Yet  in 1950's   Serre  \cite{Ser1} proved
that  some truly  non-commutative  rings $B$   satisfy  an analog of the 
fundamental  duality between polynomial rings and varieties:
\begin{equation}\label{eq1}
\mathbf{Coh}~(X)\cong \mathbf{Mod}~(B)~/~\mathbf{Tors}, 
\end{equation}
where  $X$ is a projective variety, $\mathbf{Coh}~(X)$ a category of  the quasi-coherent
sheaves on $X$, $\mathbf{Mod}~(B)$ a category of 
the finitely generated graded modules  over $B$ and  $\mathbf{Tors}$ is the torsion
[Stafford \& van den Bergh 2001]  \cite[p. 172]{StaVdb1}.
The rings $B$ satisfying (\ref{eq1}) are called  twisted  homogeneous coordinate rings,
see  Section 2 for an exact definition.   For the sake of clarity,  consider the following analogy. 
If $X$ is  a compact Hausdorff  space and $C(X)$ the commutative algebra of  
continuous functions from $X$ to ${\Bbb C}$ then by the Gelfand Duality the topology of $X$ is determined   by
the algebra $C(X)$.   In terms of the K-theory this  can be written as  
$K_0^{top}(X)\cong K_0^{alg}(C(X))$.   Taking  the two-by-two matrices 
with entries in $C(X)$,  one  gets  an  algebra 
$C(X)\otimes M_2({\Bbb C})$;  in view of stability of the K-theory
under  tensor products, 
it holds     
$K_0^{top}(X)\cong K_0^{alg}(C(X))\cong K_0^{alg}(C(X)\otimes M_2({\Bbb C}))$.  
In other words,  the topology of $X$ is defined by 
algebra $C(X)\otimes M_2({\Bbb C})$,  which is no longer a commutative algebra.  
In algebraic geometry,  one replaces $X$
by a  projective variety, $C(X)$ by its  coordinate ring,
$C(X)\otimes M_2({\Bbb C})$ by a twisted  homogeneous  
coordinate ring  of $X$ and  $K^{top}(X)$ by a category of the 
quasi-coherent sheaves on $X$.
The simplest concrete example of  $B$ is as follows. 
\begin{example}\label{ex1}
{\bf  (\cite[p. 173]{StaVdb1})}
Let $k$ be a field and $U_{\infty}(k)$ the algebra of polynomials
over $k$ in two non-commuting variables $x_1$ and $x_2$,  and a quadratic relation
$x_1x_2-x_2x_1-x_1^2=0$;  let ${\Bbb P}^1(k)$ be the projective
line over $k$.  Then $B=U_{\infty}(k)$ and $X={\Bbb P}^1(k)$
satisfy  (\ref{eq1}).  
\end{example}
In general,  there exists a canonical non-commutative ring $B$,  attached to the
projective variety $X$ and an automorphism  $\alpha: X\to X$;
we refer the reader to [Stafford \& van den Bergh 2001]  \cite[pp.180-182]{StaVdb1}.
To give an idea,  let $X=Spec~(R)$ for a commutative graded
ring $R$.  One considers the ring $B:=R[t,t^{-1}; \alpha]$ of skew
Laurent polynomials defined by the commutation relation
\begin{equation}\label{eq2}
b^{\alpha}t=tb,
\end{equation}
for all $b\in R$, where  $b^{\alpha}\in R$ is the image of $b$ under automorphism
$\alpha$;   then $B$ satisfies equation (\ref{eq1}),  see lemmas \ref{lm2} and 
\ref{lm3}.   The ring $B$ is non-commutative,  unless $\alpha$
is the  trivial automorphism of $X$.  
\begin{example}\label{ex2}
The  ring  $B=U_{\infty}(k)$
in  Example \ref{ex1}   corresponds to the automorphism $\alpha(u)=u+1$ 
of the projective line ${\Bbb P}^1(k)$.   Indeed,  $u=x_2x_1^{-1}=x_1^{-1}x_2$
and,  therefore,  $\alpha$ maps $x_2$ to $x_1+x_2$;  if  one substitutes  in (\ref{eq2})
$t=x_1,  b=x_2$ and $b^{\alpha}=x_1+x_2$,   then  one  gets the defining 
relation  $x_1x_2-x_2x_1-x_1^2=0$ for  the algebra    $U_{\infty}(k)$.   
\end{example}
In what follows,  we  consider  infinite-dimensional  representations
of $B$  by bounded linear operators  on a Hilbert space $\mathcal{H}$. 
The idea goes back to Sklyanin,  who asked about such representations
for the twisted homogeneous coordinate ring of an elliptic curve,   
see remark in brackets to the last
paragraph of [Sklyanin 1982]  \cite[Section 3]{Skl1}.

Let $\mathcal{H}$ be a Hilbert space and   $\mathcal{B}(\mathcal{H})$ the algebra of 
all  bounded linear  operators on  $\mathcal{H}$.
For a  ring of skew Laurent polynomials $R[t, t^{-1};  \alpha]$ described by
formula (\ref{eq2}),  we shall consider a homomorphism 
\begin{equation}\label{eq2bis}
\rho: R[t, t^{-1};  \alpha]\longrightarrow \mathcal{B}(\mathcal{H}). 
\end{equation}
Recall  that algebra $\mathcal{B}(\mathcal{H})$ is endowed  with a $\ast$-involution;
such an  involution  is the adjoint with respect to   the scalar product on 
the Hilbert space $\mathcal{H}$. 
\begin{definition}
We shall call representation (\ref{eq2bis})  $\ast$-coherent if:

\medskip
(i)  $\rho(t)$ and $\rho(t^{-1})$ are unitary operators,  such that
$\rho^*(t)=\rho(t^{-1})$; 

\smallskip
(ii) for all $b\in R$ it holds $(\rho^*(b))^{\alpha(\rho)}=\rho^*(b^{\alpha})$, 
where $\alpha(\rho)$ is an automorphism of  $\rho(R)$  induced by $\alpha$. 
\end{definition}
\begin{example}\label{ex2bis}
The  ring  $U_{\infty}(k)$  in examples \ref{ex1} and  \ref{ex2}  has no
$\ast$-coherent representations.
 Indeed,  involution  acts  on the generators of $U_{\infty}(k)$
 by formula $x_1^*=x_2$;  the latter does not  preserve the defining relation 
 $x_1x_2-x_2x_1-x_1^2=0$.
 \end{example}
Whenever  $B=R[t, t^{-1};  \alpha]$  admits a $\ast$-coherent representation,
$\rho(B)$ is a $\ast$-algebra;  the norm-closure of  $\rho(B)$  yields
a   $C^*$-algebra \cite[Chapter 3]{N}.  
 We shall refer to such  as   the  {\it Serre $C^*$-algebra} of the projective variety $X$;
it  will be denoted  by  $\mathcal{A}_X$.   In this case  $\alpha$ 
 induces  the Frobenius map,   see  remark \ref{rm3}.
\begin{example}\label{ex3bis}
When $X=\mathcal{E}({\Bbb C})$ a non-singular elliptic curve over the complex numbers,
 then  $B=R[t, t^{-1};  \alpha]$  is the Sklyanin algebra 
 \cite{Skl1}.   There exists a $\ast$-coherent representation
of $B$;  the resulting Serre $C^*$-algebra $\mathcal{A}_X\cong \mathcal{A}_{\theta}$,
where $\mathcal{A}_{\theta}$  is  the  noncommutative torus   \cite[Section 1.1]{N}. 
 \end{example}
Recall  that  if   $B$ is commutative,  then $X\cong \mathbf{Spec}~(B)$,   
where $\mathbf{Spec}~(B)$ is the set of all  prime  ideals of $B$.  
An objective of our  paper is to find a similar formula (if any)  
 for the Serre $C^*$-algebra $\mathcal{A}_X$.

 To get the  formula,  consider   continuous homomorphism 
 $\alpha: G\to Aut~(\mathcal{A})$,   where $G$ is a locally compact
group,  $\mathcal{A}$  is a $C^*$-algebra and $Aut~(\mathcal{A})$ the group of automorphisms of  $\mathcal{A}$. 
The triple $(\mathcal{A}, G, \alpha)$   defines a crossed product  $C^*$-algebra   denoted by 
$\mathcal{A}\rtimes_{\alpha}G$;   we refer the reader to   [Williams  2007]  \cite[pp. 47-54]{W}
for the details.   
Let $G={\Bbb Z}$ and let $\hat{\Bbb Z}\cong S^1$ be
its Pontryagin dual.  We shall write $\mathbf{Irred} ~\mathcal{A}$ for the space of all irreducible representations of  
$\mathcal{A}$ endowed with the hull-kernel (Jacobson) topology.  
Our main result can be stated as follows. 
\begin{theorem}\label{thm1}
The spaces  $\mathbf{Irred} ~(\mathcal{A}_X\rtimes_{\hat\alpha} \hat {\Bbb Z})$
and $X$ are homeomorphic, where $\hat\alpha$ is an automorphism of the Serre 
$C^*$-algebra   $\mathcal{A}_X$.
\end{theorem}
\begin{remark}
Since  $\mathcal{A}_X$ is usually  a  simple $C^*$-algebra, its spectrum
 $\mathbf{Spec}~(\mathcal{A}_X)$ is trivial. 
Thus   $X\cong  \mathbf{Irred} ~(\mathcal{A}_X\rtimes_{\hat\alpha} \hat {\Bbb Z})$  can be viewed as  a 
non-commutative analog  of the well-known isomorphism  $X\cong \mathbf{Spec}~(B)$.
\end{remark}
The article is organized as follows. The preliminary facts are reviewed 
in Section 2. Theorem \ref{thm1}  is proved in Section 3.  We consider
an application of theorem \ref{thm1}  in Section 4.

\section{Twisted homogeneous coordinate rings}
Let $X$ be a projective scheme over a field $k$, and let $\mathcal{L}$ 
be the invertible sheaf $\mathcal{O}_X(1)$ of linear forms on $X$.  Recall
that the homogeneous coordinate ring of $X$ is a graded $k$-algebra, 
which is isomorphic to the algebra
\begin{equation}\label{eq5}
B(X, \mathcal{L})=\bigoplus_{n\ge 0} H^0(X, ~\mathcal{L}^{\otimes n}). 
\end{equation}
Denote by $\mathbf{Coh}$ the category of quasi-coherent sheaves on a scheme $X$
and by $\mathbf{Mod}$ the category of graded left modules over a graded ring $B$.  
If $M=\oplus M_n$ and $M_n=0$ for $n>>0$, then the graded module
$M$ is called {\it right bounded}.  The  direct limit  $M=\lim M_{\alpha}$
is called a {\it torsion}, if each $M_{\alpha}$ is a right bounded graded
module. Denote by $\mathbf{Tors}$ the full subcategory of $\mathbf{Mod}$ of the torsion
modules.  The following result is basic about the graded ring $B=B(X, \mathcal{L})$.   
\begin{lemma}\label{lm1}
{\bf ([Serre 1955]  \cite{Ser1})}
\quad $\mathbf{Mod}~(B) ~/~\mathbf{Tors} \cong \mathbf{Coh}~(X).$
\end{lemma}

\bigskip
Let $\alpha$ be an automorphism of $X$.  The pullback of sheaf $\mathcal{L}$ 
along $\alpha$ will be denoted by $\mathcal{L}^{\alpha}$,  i.e. 
$\mathcal{L}^{\alpha}(U):= \mathcal{L}(\alpha U)$ for every $U\subset X$. 
We shall set
\begin{equation}\label{eq7}
B(X, \mathcal{L}, \alpha)=\bigoplus_{n\ge 0} H^0(X, ~\mathcal{L}\otimes \mathcal{L}^{\alpha}\otimes\dots
\otimes  \mathcal{L}^{\alpha^{ n}}). 
\end{equation}

 The multiplication of sections is defined by the rule
 \begin{equation}\label{eq8}
 ab=a\otimes b^{\alpha^m},
 \end{equation}
 whenever $a\in B_m$ and $b\in B_n$.

 Given a pair $(X,\alpha)$ consisting of a Noetherian scheme $X$ and 
 an automorphism $\alpha$ of $X$,  an invertible sheaf $\mathcal{L}$ on $X$
 is called {\it $\alpha$-ample}, if for every coherent sheaf $\mathcal{F}$ on $X$,
 the cohomology group $H^q(X, ~\mathcal{L}\otimes \mathcal{L}^{\alpha}\otimes\dots
\otimes  \mathcal{L}^{\alpha^{ n-1}}\otimes \mathcal{F})$  vanishes for $q>0$ and
$n>>0$.  Notice,  that if $\alpha$ is trivial,  this definition is equivalent to the
usual definition of ample invertible sheaf [Serre 1955]  \cite{Ser1}.    
A  non-commutative generalization of the Serre theorem is as follows.
\begin{lemma}\label{lm2}
{\bf ([Artin \& van den Bergh  1990]  \cite{ArtVdb1})}
Let $\alpha: X\to X$ be an automorphism of a projective scheme $X$
over $k$  and let $\mathcal{L}$ be a $\alpha$-ample invertible sheaf on $X$. If
$B(X, \mathcal{L}, \alpha)$  is  the ring (\ref{eq7}),    then
\begin{equation}\label{eq9}
\mathbf{Mod}~(B(X, \mathcal{L}, \alpha)) ~/~\mathbf{Tors} \cong \mathbf{Coh}~(X).  
\end{equation}
\end{lemma}
\begin{remark}\label{rm1.5}
\textnormal{
The necessary and sufficient conditions for an invertible sheaf to be 
$\alpha$-ample have been established by  [Keeler 2000] \cite{Kee1}. 
}
\end{remark}

\section{Proof of theorem \ref{thm1}}
We shall split the proof in a series of lemmas, starting with the following
\begin{lemma}\label{lm3}
 $B(X, \mathcal{L}, \alpha)\cong R[t,t^{-1}; \alpha]$,  where $X=\mathbf{Spec}~(R)$
 of the commutative ring $R$.
\end{lemma}
\begin{proof}
Let us  write the twisted homogeneous coordinate ring  $B(X, \mathcal{L}, \alpha)$ of projective variety
$X$ in the following  form:
\begin{equation}\label{eq10}
B(X, \mathcal{L}, \alpha)=\bigoplus_{n\ge 0} H^0(X, \mathfrak{B}_n),
\end{equation}
where $\mathfrak{B}_n=\mathcal{L}\otimes \mathcal{L}^{\alpha}\otimes\dots
\otimes  \mathcal{L}^{\alpha^{ n}}$ and  $H^0(X, \mathfrak{B}_n)$ is the zero
sheaf cohomology of  $X$,  i.e. the space of sections $\Gamma(X, \mathfrak{B}_n)$;
compare with formula (3.5) of [Artin \& van den Bergh  1990]  \cite{ArtVdb1}.

If one denotes by $\mathcal{O}$  the structure sheaf of $X$, then 
\begin{equation}\label{eq10bis}
\mathfrak{B}_n=\mathcal{O}t^n
\end{equation}
can be interpreted  as a free left  $\mathcal{O}$-module of rank one with basis $\{t^n\}$ 
 [Artin \& van den Bergh  1990]  \cite{ArtVdb1},  p. 252.

Recall, that spaces  $B_i=H^0(X, \mathfrak{B}_i)$ have  been endowed with the multiplication
rule  (\ref{eq8}) between the   sections  $a\in B_m$ and $b\in B_n$;
such a rule translates into the formula:
\begin{equation}\label{eq11}
at^mbt^n=ab^{\alpha^m}t^{m+n}. 
\end{equation}

One can eliminate $a$ and $t^n$ in   the both sides of (\ref{eq11});  
this operation gives us the following equation: 
\begin{equation}\label{eq12}
t^mb=b^{\alpha^m}t^m. 
\end{equation}

First notice, that our ring  $B(X, \mathcal{L}, \alpha)$ contains a commutative
subring $R$, such that $\mathbf{Spec}~(R)=X$. 
Indeed, let $m=0$ in formula (\ref{eq12});  then $b=b^{Id}$ and, thus, $\alpha=Id$.
We conclude therefore, that $R=B_0$ is a commutative subring of  
$B(X, \mathcal{L}, \alpha)$,  and $\mathbf{Spec}~(R)=X$.

Let us show that equations  (\ref{eq2}) and  (\ref{eq12}) are equivalent. 
First, let us show that  (\ref{eq2}) implies  (\ref{eq12}).   Indeed,  equation 
(\ref{eq2}) can be written as $b^{\alpha}=tbt^{-1}$.   Then: 
\begin{equation}\label{eq12bis}
\left\{
\begin{array}{ccc}
b^{\alpha^2} &=& tb^{\alpha}t^{-1}= t^2 bt^{-2},\\
b^{\alpha^3} &=& tb^{\alpha^2}t^{-1}= t^3 b t^{-3},\\
&\vdots&\\
b^{\alpha^m} &=& tb^{\alpha^{m-1}}t^{-1}= t^m b t^{-m}. 
\end{array}
\right.
\end{equation}
The last line of (\ref{eq12bis}) is  equivalent to equation (\ref{eq12}). 
The converse is evident;  one sets $m=1$ in (\ref{eq12}) and obtains
equation (\ref{eq2}).  Thus,    (\ref{eq2}) and  (\ref{eq12}) are equivalent
equations.

It is easy now to establish an isomorphism  $B(X, \mathcal{L}, \alpha)\cong R[t,t^{-1}; \alpha]$.
For that,  take  $b\in R\subset B(X, \mathcal{L}, \alpha)$;  then $B(X, \mathcal{L}, \alpha)$ 
coincides  with the ring of the skew  Laurent polynomials  $R[t,t^{-1}; \alpha]$,
since the  commutation relation (\ref{eq2})  is equivalent to equation (\ref{eq12}).
 Lemma \ref{lm3} follows. 
 \end{proof}
\begin{remark}
Lemma \ref{lm3} is proved for an affine variety.  For such varieties the theory 
of twisted homogeneous coordinate rings can be  extended  by replacing an ${\Bbb N}$-grading
by  the ${\Bbb Z}$-grading. 
\end{remark}

 \bigskip
\begin{lemma}\label{lm4}
$\mathcal{A}_X\cong C(X)\rtimes_{\alpha} {\Bbb Z}$,  where $C(X)$ is the
$C^*$-algebra of all continuous complex-valued functions on $X$
and $\alpha$  is a $\ast$-coherent  automorphism of  $X$.   
\end{lemma}
\begin{proof}
By definition of the Serre algebra $\mathcal{A}_X$,   the ring of skew Laurent 
polynomials $R[t, t^{-1}; \alpha]$  is  dense in $\mathcal{A}_X$;  roughly
speaking, one has to show that this property defines a crossed product
structure on $\mathcal{A}_X$.  We shall proceed in the following steps.

\medskip
(i) Recall  that  $R[t, t^{-1}; \alpha]$ consists of the finite sums
\begin{equation}\label{eq13}
\sum  b_k t^k,  \qquad b_k\in R,
\end{equation}

 subject to the commutation relation
\begin{equation}\label{eq13bis}
b_k^{\alpha}t=tb_k. 
\end{equation}
Thanks to a $\ast$-coherent representation,
there is also an involution on $R[t, t^{-1}; \alpha]$, subject to the
following rules:
\begin{equation}\label{eq15}
\left\{
\begin{array}{ccccc}
 &(i)&   t^* &=& t^{-1},\\
 &(ii)&  (b_k^*)^{\alpha} &=& (b_k^{\alpha})^*.
\end{array}
\right.
\end{equation}

\bigskip
(ii)  Following [Williams  2007]  \cite[p.47]{W},  we shall consider the set 
$C_c({\Bbb Z}, R)$ of continuous functions from 
${\Bbb Z}$ to $R$ having a compact support;  
then the formal sums (\ref{eq13}) can be viewed as 
elements of  $C_c({\Bbb Z}, R)$ via  the identification 
\begin{equation}\label{eq15bis}
k\longmapsto b_k.
\end{equation}

It can be verified  that multiplication operation of the formal sums (\ref{eq13})
translates into a convolution product of functions $f,g\in C_c({\Bbb Z}, R)$
given by the formula:
\begin{equation}\label{eq15bisbis}
(f g)(k) = \sum_{l\in {\Bbb Z}} f(l) t^l g(k-l) t^{-l}, 
\end{equation}
while involution (\ref{eq15}) translates into an involution on $C_c({\Bbb Z}, R)$
 given by the formula:
\begin{equation}\label{eq16}
f^*(k) = t^k f^*(-k) t^{-k}.
\end{equation}

The multiplication given by convolution product (\ref{eq15bisbis}) 
and involution (\ref{eq16}) turn $C_c({\Bbb Z}, R)$ into a $\ast$-algebra,
which is isomorphic to  the algebra $R[t, t^{-1}; \alpha]$.

\bigskip
(iii)  There exists the standard construction of a norm on  
$C_c({\Bbb Z}, R)$;  we omit it here referring the reader 
to [Williams  2007]  \cite{W}, Section 2.3.  The completion of $C_c({\Bbb Z}, R)$
in that norm defines a crossed product $C^*$-algebra 
$R\rtimes_{\alpha}{\Bbb Z}$ [Williams  2007]  \cite{W}, Lemma 2.27.

\bigskip
(iv)  Since $R$ is a commutative $C^*$-algebra
and $X=\mathbf{Spec}~(R)$,  one concludes that $R\cong C(X)$.
Thus, one obtains $\mathcal{A}_X=C(X)\rtimes_{\alpha}{\Bbb Z}$.

Lemma \ref{lm4} follows. 
\end{proof}
\begin{remark}\label{rm2}
It is easy to see,   that (\ref{eq13bis}) and (\ref{eq15}i) imply (\ref{eq15}ii);
in other words,  if involution does not commute with automorphism $\alpha$,
representation $\rho$  cannot  be  unitary,  i.e. $\rho^*(t)\ne\rho(t^{-1})$.    
 \end{remark}

\begin{lemma}\label{lm5}
There exists  $\hat\alpha\in Aut~(\mathcal{A}_X)$, such that: 
\begin{equation}\label{eq19}
X\cong \mathbf{Irred}~(\mathcal{A}_X\rtimes_{\hat\alpha} \hat {\Bbb Z}).
\end{equation}
\end{lemma}
\begin{proof}
 Formula (\ref{eq19}) is an implication of the Takai duality 
for the crossed products  [Williams  2007]  \cite[Section 7.1]{W}; let us briefly 
review this construction.  

Let $(A,G,\alpha)$ be a $C^*$-dynamical system with $G$ locally compact
abelian group; let $\hat G$ be the dual of $G$. For each $\gamma\in \hat G$,
one can define a map $\hat a_{\gamma}: C_c(G,A)\to C_c(G,A)$
given by the formula:
\begin{equation}\label{eq20}
\hat a_{\gamma}(f)(s)=\bar\gamma(s)f(s), \qquad\forall s\in G.
 \end{equation}
In fact, $\hat a_{\gamma}$ is a $\ast$-homomorphism, since it
respects the convolution product and involution on $C_c(G,A)$
[Williams  2007]  \cite{W}.  Because the crossed product $A\rtimes_{\alpha}G$
is the closure of $C_c(G,A)$, one gets an extension of $\hat a_{\gamma}$
to an element of $Aut~(A\rtimes_{\alpha}G)$ and, therefore, a 
homomorphism:
\begin{equation}\label{eq21}
\hat\alpha: \hat G\to Aut~(A\rtimes_{\alpha}G).
 \end{equation}
Recall that the  Takai duality asserts that:
\begin{equation}\label{eq22}
(A\rtimes_{\alpha} G)\rtimes_{\hat\alpha}\hat G\cong A\otimes \mathcal{K}(L^2(G)),
 \end{equation}
where $\mathcal{K}(L^2(G))$ is the algebra of compact operators on the
Hilbert space $L^2(G)$. 

Let us substitute $A=C_0(X)$   and   $G={\Bbb Z}$ in (\ref{eq22}); 
one gets the following isomorphism:
\begin{equation}\label{eq23}
(C_0(X)\rtimes_{\alpha} {\Bbb Z})\rtimes_{\hat\alpha}\hat {\Bbb Z}\cong C_0(X)\otimes 
\mathcal{K}(L^2({\Bbb Z})).
 \end{equation}

Lemma \ref{lm4} says that $C_0(X)\rtimes_{\alpha}{\Bbb Z}\cong \mathcal{A}_X$;
 therefore  one arrives at the following isomorphism: 
\begin{equation}\label{eq24}
\mathcal{A}_X \rtimes_{\hat\alpha} \hat {\Bbb Z}\cong C_0(X)\otimes 
\mathcal{K}(L^2({\Bbb Z})).
 \end{equation}

Consider the set of all irreducible representations of the $C^*$-algebras
in (\ref{eq24});  then one gets the following equality of representations:
\begin{equation}\label{eq25}
\mathbf{Irred}~(\mathcal{A}_X \rtimes_{\hat\alpha} \hat {\Bbb Z}) =\mathbf{Irred}~(C_0(X)\otimes 
\mathcal{K}(L^2({\Bbb Z}))).
 \end{equation}

Let $\pi$ be a representation of the tensor product 
$C_0(X)\otimes \mathcal{K}(L^2({\Bbb Z}))$ on the Hilbert
space $\mathcal{H}\otimes L^2({\Bbb Z})$;  then $\pi=\varphi\otimes\psi$,
where $\varphi: C_0(X)\to \mathcal{B}(\mathcal{H})$ and $\psi: \mathcal{K}\to \mathcal{B}(L^2({\Bbb Z}))$. 
It is known, that the only irreducible representation of the algebra of 
compact operators is the identity representation. Thus, one gets:
\begin{eqnarray}\label{eq26}
\mathbf{Irred}~(C_0(X)\otimes \mathcal{K}(L^2({\Bbb Z}))) &=&
\mathbf{Irred}~(C_0(X))\otimes \{pt\}=\\
       &=& \mathbf{Irred}(C_0(X)).\nonumber
 \end{eqnarray}

Further, the $C^*$-algebra $C_0(X)$ is commutative,  hence the 
following equations are true:
\begin{equation}\label{eq27}
\mathbf{Irred}~(C_0(X))=\mathbf{Spec}~(C_0(X))=X.
 \end{equation}
 
 Putting together (\ref{eq25}) -- (\ref{eq27}),   one obtains:
\begin{equation}\label{eq28}
\mathbf{Irred}~(\mathcal{A}_X \rtimes_{\hat\alpha} \hat {\Bbb Z})  \cong    X.
 \end{equation}

The conclusion of lemma \ref{lm5} follows from (\ref{eq28}). 
\end{proof}

\bigskip
Theorem \ref{thm1} follows from lemma \ref{lm5}.

\section{Example}
Recall that noncommutative torus  is  the  universal $C^*$-algebra   $\mathcal{A}_{\theta}$ generated by 
   unitaries  $u$ and $v$ satisfying the commutation  relation $vu=e^{2\pi i\theta}uv$, where   $\theta$ is an irrational constant
 \cite[Section 1.1]{N}. 
The K-theory of $\mathcal{A}_{\theta}$ is  Bott periodic with  $K_0(\mathcal{A}_{\theta})=K_1(\mathcal{A}_{\theta})\cong {\Bbb Z}^2$;
the range of trace on  projections of $\mathcal{A}_{\theta}\otimes \mathcal{K}$ is a subset
$\Lambda={\Bbb Z}+{\Bbb Z}\theta$  of the real line called a 
pseudo-lattice.   
The torus $\mathcal{A}_{\theta}$ is said to have   real multiplication,
if $\theta$ is a quadratic irrationality;  in this case the endomorphism ring of pseudo-lattice
$\Lambda$ is  bigger than ${\Bbb Z}$ -- hence the name.  
The noncommutative torus  with  real multiplication will be written as $\mathcal{A}_{RM}$.  
We denote by  $(\overline{a_1,a_2,\dots, a_n})$ be the minimal period 
of continued fraction of the irrational quadratic $\theta$. 
Consider the  matrix
\begin{equation}
A:= \left(
\begin{matrix}
a_1 & 1\cr 1 & 0
\end{matrix}\right) 
\left(
\begin{matrix}a_2 & 1\cr 1 & 0
\end{matrix}
\right)\dots
\left(
\begin{matrix}
a_n & 1\cr 1 & 0
\end{matrix}
\right);
\end{equation}
the matrix is an  invariant of  torus  $\mathcal{A}_{RM}$. 
A specialization of theorem \ref{thm1} to $\mathcal{A}_{RM}$
gives us the following result.
\begin{proposition}\label{ex3}
\quad$\mathbf{Irred}~(\mathcal{A}_{RM}\rtimes_{\hat\alpha} \hat {\Bbb Z})\cong \mathcal{E}(K)$,
where $\mathcal{E}(K)$ is a non-singular elliptic curve defined over an algebraic 
number field $K$. 
\end{proposition}
\begin{proof}
We shall view the crossed product $\mathcal{A}_{RM}\rtimes_{\hat\alpha} \hat {\Bbb Z}$
as a $C^*$-dynamical system $(\mathcal{A}_{RM}, \hat {\Bbb Z}, \hat\alpha)$,  
see [Williams  2007]  \cite{W} for the details. 
Recall that the irreducible representations of 
$C^*$-dynamical system  $(\mathcal{A}_{RM}, \hat {\Bbb Z}, \hat\alpha)$
 are in the one-to-one correspondence with the minimal sets of the 
dynamical system (i.e. closed $\hat\alpha$-invariant sub-$C^*$-algebras
of $\mathcal{A}_{RM}$ not containing a smaller object with the same property).

To calculate the minimal sets of $(\mathcal{A}_{RM}, \hat {\Bbb Z}, \hat\alpha)$,
let $\theta$ be quadratic irrationality such that $\mathcal{A}_{RM}\cong \mathcal{A}_{\theta}$.
It is known that every non-trivial sub-$C^*$-algebra  of $\mathcal{A}_{\theta}$ has
the form $\mathcal{A}_{n\theta}$ for some positive integer $n$
 [Rieffel 1981] \cite[p. 419]{Rie2}.   It is easy to deduce that the
{\it maximal}  proper sub-$C^*$-algebra of $\mathcal{A}_{\theta}$ has the form 
$\mathcal{A}_{p\theta}$, where $p$ is   a prime number.  
(Indeed, each composite $n=n_1n_2$ cannot be maximal since
$\mathcal{A}_{n_1n_2\theta}\subset \mathcal{A}_{n_1\theta}\subset \mathcal{A}_{\theta}$
or $\mathcal{A}_{n_1n_2\theta}\subset \mathcal{A}_{n_2\theta}\subset \mathcal{A}_{\theta}$,
where all inclusions are strict.)

We claim that  $(\mathcal{A}_{p\theta}, \hat {\Bbb Z}, \hat\alpha^{\pi(p)})$
is the minimal $C^*$-dynamical system,  where $\pi(p)$ is certain
power of the automorphism $\hat\alpha$.   Indeed,  the automorphism 
$\hat\alpha$ of $\mathcal{A}_{\theta}$  corresponds to multiplication by
the fundamental unit, $\varepsilon$, of   pseudo-lattice  $\Lambda=
{\Bbb Z}+\theta {\Bbb Z}$.   It is known that certain power, $\pi(p)$,
of $\varepsilon$  coincides with the fundamental unit of pseudo-lattice
${\Bbb Z}+(p\theta){\Bbb Z}$, see e.g. [Hasse 1950] \cite[p.298]{H}.  
Thus  one gets the minimal $C^*$-dynamical system 
$(\mathcal{A}_{p\theta}, \hat {\Bbb Z}, \hat\alpha^{\pi(p)})$,  which is defined 
on the sub-$C^*$-algebra $\mathcal{A}_{p\theta}$ of  $\mathcal{A}_{\theta}$.  
Therefore we have an isomorphism    
\begin{equation}\label{eq29}
\mathbf{Irred} ~(\mathcal{A}_{RM}\rtimes_{\hat\alpha} \hat {\Bbb Z})\cong
\bigcup_{p\in \mathcal{P}} \mathbf{Irred} ~( \mathcal{A}_{p\theta}\rtimes_{\hat\alpha^{\pi(p)}} \hat {\Bbb Z}),
 \end{equation}
where $\mathcal{P}$ is the set of all (but a finite number) of primes.

To simplify the RHS of (\ref{eq29}), 
recall that matrix form of the fundamental unit $\varepsilon$ of
pseudo-lattice $\Lambda$ coincides with the matrix $A$, see above. 
For each prime $p\in \mathcal{P}$ consider the matrix 
\begin{equation}\label{eq30}
L_p=\left(
\begin{matrix}tr~(A^{\pi(p)})-p & p\cr tr~(A^{\pi(p)})-p-1 & p
\end{matrix}\right),
 \end{equation}
where $tr$ is the trace of matrix.  Let us show, that 
\begin{equation}\label{eq31}
\mathcal{A}_{p\theta}\rtimes_{\hat\alpha^{\pi(p)}} \hat {\Bbb Z}\cong
\mathcal{A}_{\theta}\rtimes_{L_p} \hat {\Bbb Z},
 \end{equation}
where $L_p$ is an endomorphism of $\mathcal{A}_{\theta}$ (of degree $p$)
induced by matrix $L_p$.
Indeed,  because $deg~(L_p)=p$ the endomorphism $L_p$ maps
pseudo-lattice $\Lambda={\Bbb Z}+\theta {\Bbb Z}$ to a sub-lattice
of index $p$;  any such can be written in the form 
$\Lambda_p={\Bbb Z}+(p\theta) {\Bbb Z}$, see e.g. 
[Borevich \& Shafarevich  1966]  \cite[p.131]{BS}.
Notice that pseudo-lattice $\Lambda_p$ corresponds
to the sub-$C^*$-algebra $\mathcal{A}_{p\theta}$  of algebra $\mathcal{A}_{\theta}$  
and $L_p$ induces a shift automorphism of $\mathcal{A}_{p\theta}$ 
 [Cuntz 1977]  \cite[Section 2.1]{Cun1}.
It is not hard to see, that the shift automorphism coincides
with $\hat\alpha^{\pi(p)}$.   Indeed, it is verified directly
that $tr~(\hat\alpha^{\pi(p)})=tr~(A^{\pi(p)})=tr~(L_p)$;
thus one gets a bijection between powers of $\hat\alpha^{\pi(p)}$
and such of $L_p$.   But $\hat\alpha^{\pi(p)}$ corresponds to
the fundamental unit of pseudo-lattice $\Lambda_p$;  therefore 
the shift automorphism induced by $L_p$
must coincide with $\hat\alpha^{\pi(p)}$.   The isomorphism (\ref{eq31})
is proved.

Therefore (\ref{eq29}) can be simplified to the form
\begin{equation}\label{eq32}
\mathbf{Irred} ~(\mathcal{A}_{RM}\rtimes_{\hat\alpha} \hat {\Bbb Z})\cong
\bigcup_{p\in \mathcal{P}} \mathbf{Irred} ~( \mathcal{A}_{RM}\rtimes_{L_p} \hat {\Bbb Z}).  
 \end{equation}

To calculate  irreducible representations  of the crossed product $C^*$-algebra 
 $\mathcal{A}_{RM}\rtimes_{L_p} \hat {\Bbb Z}$ at the RHS of  (\ref{eq32}),
 recall that such are in a one-to-one correspondence with 
 the set of  invariant measures on a subshift of finite type given by the positive integer 
 matrix (\ref{eq30})  [Bowen \& Franks 1977]  \cite{BowFra1};  the  measures make an  abelian group under the 
  addition operation.   Such a  group is  isomorphic to  ${\Bbb Z}^2~/~(I-L_p) {\Bbb Z}^2$,
where $I$ is the identity matrix, see  [Bowen \& Franks 1977]  \cite[Theorem 2.2]{BowFra1}.

Therefore (\ref{eq32}) can be written in  the form
\begin{equation}\label{eq33}
\mathbf{Irred} ~(\mathcal{A}_{RM}\rtimes_{\hat\alpha} \hat {\Bbb Z})\cong
\bigcup_{p\in \mathcal{P}} {{\Bbb Z}^2\over (I-L_p) {\Bbb Z}^2}.
\end{equation}

Let $\mathcal{E}(K)$ be a non-singular elliptic curve defined over 
the algebraic number field $K$;    let $\mathcal{E}({\Bbb F}_p)$ 
be the reduction of $\mathcal{E}(K)$ modulo prime ideal over a 
``good'' prime number $p$.  Recall that $|\mathcal{E}({\Bbb F}_p)|=
det~(I-Fr_p)$,  where  $Fr_p$  is an integer two-by-two matrix
corresponding to the action of Frobenius endomorphism on the
$\ell$-adic cohomology of $\mathcal{E}(K)$,  see e.g. 
[Tate  1974]  \cite[p.187]{Tat1}.

Since $|{\Bbb Z}^2 / (I-L_p){\Bbb Z}^2|=det~(I-L_p)$,
one can identify $Fr_p$ and $L_p$  and,  therefore,  
one obtains an isomorphism $\mathcal{E}({\Bbb F}_p)\cong 
{\Bbb Z}^2 / (I-L_p){\Bbb Z}^2$.  
Thus  (\ref{eq33}) can be written in  the form
\begin{equation}\label{eq34}
\mathbf{Irred} ~(\mathcal{A}_{RM}\rtimes_{\hat\alpha} \hat {\Bbb Z})\cong
\bigcup_{p\in \mathcal{P}} \mathcal{E}({\Bbb F}_p).
\end{equation}

Finally,   consider an arithmetic scheme, $X$,  corresponding to $\mathcal{E}(K)$;
the latter fibers over ${\Bbb Z}$ [Silverman 1994]  \cite[Example 4.2.2]{S}.
   It can be immediately seen,  that the RHS of (\ref{eq34})
coincides with the scheme $X$, where  the regular fiber over $p$ corresponds
to $\mathcal{E}({\Bbb F}_p)$ {\it ibid.}  This argument finishes the proof of 
proposition \ref{ex3}.
\end{proof}
\begin{remark}\label{rm3}
The Frobenius endomorphism $Fr_p=L_p$
is induced by map $\alpha$,   since $\alpha$ defines $\hat\alpha$ 
(i.e. matrix $A$) and the latter is linked to  $L_p$  via  formula (\ref{eq30}).
\end{remark}


\bibliographystyle{amsplain}

\begin{thebibliography}{99}



 \bibitem{ArtVdb1}
M.~Artin and M. ~van den Bergh,  
\textit{Twisted homogeneous coordinate rings}, J. of Algebra {\bf 133} 
(1990), 249-271.  





\bibitem{BS}
Z.~I.~Borevich and I.~R.~Shafarevich,  
\textit{Number Theory}, 
Academic Press, 1966.    


\bibitem{BowFra1}
R.~Bowen and J.~Franks, 
\textit{Homology for zero-dimensional nonwandering
sets},  Ann. of Math. {\bf 106} (1977), 73-92.  


\bibitem{Cun1}
J.~Cuntz,   
\textit{Simple $C^*$-algebras generated by isometries}, 
Commun. Math. Phys. {\bf 57} (1977), 173-185. 





\bibitem{H}
H.~Hasse,   
\textit{Vorlesungen \"uber   Zahlentheorie},   Springer, 1950.


\bibitem{Kee1} 
D.~S.~Keeler, 
\textit{Criteria  for  $\sigma$-ampleness},
J. ~Amer. Math. Soc. {\bf 13} (2000),  517-532. 


\bibitem{N}
I.~V.~Nikolaev, \textit{Noncommutative Geometry},
De Gruyter Studies in Math. {\bf 66}, Berlin, 2017.


\bibitem{Rie2} 
M.~A.~Rieffel, 
\textit{$C^*$-algebras associated with irrational rotations},
Pacific J. Math. {\bf 93} (1981),  415-429.  



\bibitem{Ser1}
J.~P.~Serre,  
\textit{Fasceaux alg\'ebriques coh\'erents}, Ann. of Math. {\bf 61} (1955), 
197-278.  



\bibitem{S}
J.~H.~Silverman, 
\textit{Advanced Topics in the Arithmetic of Elliptic Curves},
GTM {\bf 151}, Springer 1994.




\bibitem{Skl1}
E.~K.~Sklyanin, 
\textit{Some algebraic structures connected to the Yang--Baxter
equations}, Functional Anal. Appl. {\bf 16} (1982), 27-34. 



\bibitem{StaVdb1}
J.~T.~Stafford and M.~van ~den ~Bergh, 
\textit{Noncommutative curves and noncommutative
surfaces}, Bull. Amer. Math. Soc. {\bf 38} (2001), 171-216. 


\bibitem{Tat1}
J.~T.~Tate, 
\textit{The arithmetic of elliptic curves}, Inventiones  Math. {\bf 23} (1974),
179-206. 


\bibitem{W}
D.~P.~Williams,  
\textit{Crossed Products of $C^*$-Algebras}, Math. Surveys and 
Monographs, Vol. {\bf 134},  Amer. Math. Soc.  2007.
   

\end{thebibliography}


\end{document}